\documentclass [12pt]{article}
\usepackage{graphicx,amssymb,amsfonts,latexsym,amsmath,amsthm,times}
\usepackage{epsfig}
\usepackage{color}
\usepackage[english]{babel}
\usepackage[a4paper]{geometry}
\def\lanbox{\hbox{$\, \vrule height 0.25cm width 0.25cm depth 0.01cm \,$}}

\numberwithin{equation}{section}

\begin{document}

\vspace*{1.4cm}

\normalsize \centerline{\Large \bf ON THE EXISTENCE OF
STATIONARY SOLUTIONS FOR} 

\medskip

\normalsize \centerline{\Large \bf CERTAIN
SYSTEMS OF INTEGRO-DIFFERENTIAL EQUATIONS}

\medskip

\centerline{\Large\bf WITH THE DOUBLE SCALE ANOMALOUS DIFFUSION}

\vspace*{1cm}

\centerline{\bf Vitali Vougalter$^{1 \ *}$,  Vitaly Volpert$^{2, 3}$,}

\vspace*{0.5cm}

\centerline{$^{1 \ *}$ Department of Mathematics, University
of Toronto}

\centerline{Toronto, Ontario, M5S 2E4, Canada}

\centerline{ e-mail: vitali@math.toronto.edu}

\medskip

\centerline{$^2$ Institute Camille Jordan, UMR 5208 CNRS,
University Lyon 1}

\centerline{ Villeurbanne, 69622, France}

\centerline{$^3$ Peoples' Friendship University of Russia, 6 Miklukho-Maklaya
St,}

\centerline{Moscow, 117198, Russia}

\centerline{e-mail: volpert@math.univ-lyon1.fr}

\medskip

%*******************************************************************
%ABSTRACT
%*******************************************************************

\vspace*{0.25cm}

\noindent {\bf Abstract:}
The work deals with establishing the solvability of a system of
integro-differential equations in the situation of the double scale anomalous
diffusion. Each equation of such system involves the sum of the two negative
Laplace operators raised to two distinct fractional powers in the space
of three dimensions. The proof of the existence of solutions is based
on a fixed point technique. We use the solvability conditions for the
non-Fredholm elliptic operators in unbounded domains.

\vspace*{0.25cm}



\vspace*{0.5cm}

\bigskip

\bigskip

%%%%%%%%%%%%%%%%%%%%%%%%%%%%%%%%%%%%%%%%%%%%%%%%%%

\setcounter{section}{1}

\centerline{\bf 1. Introduction}

\medskip

In the present article we study the existence of stationary
solutions of the following system of $N\geq 2$ integro-differential equations in
${\mathbb R}^{3}$
$$
\frac{\partial u_{m}}{\partial t}=
-D_{m}[(-\Delta)^{s_{1, m}}+ (-\Delta)^{s_{2, m}}]u_{m} +
$$
\begin{equation}
\label{h}
\int_{{\mathbb R}^{3}}K_{m}(x-y)g_{m}(u(y,t))dy + f_{m}(x), 
\end{equation}
with $\displaystyle{1\leq m\leq N, \ \frac{1}{4}<s_{1, m}<\frac{3}{4}}$ and
$\displaystyle{{s}_{1, m}<s_{2, m}<1}$, which appears in the cell
population dynamics.
The results of the work are derived in these particular ranges of the
values of the powers of the negative Laplacians, which is based on the
solvability of the linear Poisson type equations (\ref{lpm}) and the
applicability of the Sobolev inequality (\ref{sobfr}) for the fractional
Laplace operator. The solvability of the single equation analogous to
(\ref{h}) involving a single fractional Laplacian in the diffusion term
was discussed in ~\cite{VV15}.
Let us note that the space variable $x$ in our system is 
correspondent to the cell genotype, the functions
$u_{m}(x,t)$ describe the cell density distributions for various groups of
cells as functions of their genotype and time,
$$
u(x,t)=(u_{1}(x,t), u_{2}(x,t),...,u_{N}(x,t))^{T}.
$$
The right side of the system of equations (\ref{h}) describes the evolution of
cell densities by means of the cell proliferation, mutations and the cell influx
or efflux. The double scale anomalous diffusion terms with positive
coefficients $D_m$ correspond to the change of genotype due to small random
mutations, and the
integral production terms describe large mutations.
Functions $g_{m}(u)$ stand for the rates of the cell birth depending on $u$
(density dependent proliferation), and the kernels $K_{m}(x-y)$ express
the proportions of the newly born cells, which  change their genotypes from $y$ to $x$.
We assume that they depend on the distance between the genotypes.
The functions $f_{m}(x)$ denote the influxes or effluxes of cells for
different genotypes.

The fractional Laplacian describes
a particular case of the anomalous diffusion actively studied in the context of
the various applications in plasma physics and
turbulence \cite{2}, \cite{1}, surface diffusion \cite{3}, \cite{4},
semiconductors \cite{5} and so on. The anomalous
diffusion can be understood as a random process of the particle motion
characterized by the probability density distribution of the jump length.
The moments of this density distribution are finite in
the case of the normal diffusion, but this is not the case for the anomalous
diffusion. The asymptotic behavior at the infinity of the probability density
function determines the value of the power of the negative Laplacian
(see ~\cite{MK00}). Weak error for continuous time Markov chains related to
fractional in time P(I)DEs was estimated in ~\cite{KKM16}.
In the present work we consider the case of
$\displaystyle{\frac{1}{4}<s_{1, m}<\frac{3}{4}, \ s_{1, m}<s_{2, m}<1, \
1\leq m\leq N}$. The necessary conditions of the preservation
of the nonnegativity of the solutions of a system of parabolic equations in
the case of the double scale anomalous diffusion were derived in
~\cite{EV221}. In the article ~\cite{GNTUV21} the authors discuss the
simultaneous inversion for the fractional exponents in the space-time
fractional diffusion equation. Stability results for evolution equations 
involving superposition operators of mixed fractional orders were
obtained in ~\cite{KK25}. An existence theory for superposition operators
of mixed order subject to jumping nonlinearities was developed in ~\cite{DKSV24}.

Let us set here all $D_{m}=1$ and establish the existence of solutions of the
system of equations 
\begin{equation}
\label{p}
-[(-\Delta)^{s_{1, m}}+(-\Delta)^{s_{2, m}}]u_{m}+\int_{{\mathbb R}^{3}} K_{m}(x-y)g_{m}
(u(y))dy +f_{m}(x) = 0,
\end{equation}
where
$\displaystyle{\frac{1}{4}<s_{1, m}<\frac{3}{4}, \ s_{1, m}<s_{2, m}<1, \
1\leq m\leq N}$.
We discuss the case when the linear parts of the operators involved in
our system of equations does not satisfy the Fredholm property. As a consequence,  the
conventional methods of the nonlinear analysis may not be used.  Our argument is based on
the solvability conditions for the operators without the Fredholm property and
the method of contraction mappings.

Consider the equation
\begin{equation}
\label{eq1}
 -\Delta u + V(x) u - a u=f,
\end{equation}
with $u \in E= H^{2}({\mathbb R}^{d})$ and  $f \in F=
L^{2}({\mathbb R}^{d}), \ d\in {\mathbb N}$, $a$ is a constant and
the scalar potential function $V(x)$ is either trivial in the whole space
or tends to $0$ as $|x|\to \infty$. This model problem is discussed here
to illustrate the particular features of the equations containing the non-Fredholm
operators, the techniques used to solve them and the previous results.
If $a \geq 0$, the essential spectrum of the
operator $A : E \to F$ corresponding to the left side of problem
(\ref{eq1}) contains the origin. As a consequence, this operator fails to
satisfy the Fredholm property. Its image is not closed, for $d>1$
the dimension of its kernel and the codimension of its image are
not finite. The present work is devoted to the studies of the certain properties
of the operators of this kind. The elliptic problems involving the non-Fredholm
operators were considered actively in recent years.
Approaches in weighted Sobolev and H\"older spaces were developed in
~\cite{Amrouche1997}, ~\cite{Amrouche2008}, ~\cite{Bolley1993},
~\cite{Bolley2001}, ~\cite{B88}. The Schr\"odinger type operators without the
Fredholm property were treated with the methods of the spectral and the
scattering theory in ~\cite{EV22}, ~\cite{V2011}, ~\cite{VV08}, ~\cite{VV19}.
The nonlinear non-Fredholm elliptic problems were discussed in ~\cite{EV22},
~\cite{EV221}, ~\cite{EV24}, ~\cite{VV111}, ~\cite{VV14}, ~\cite{VV15}. The
significant applications to the theory of reaction-diffusion
type equations were investigated in ~\cite{DMV05}, ~\cite{DMV08}. Fredholm
structures, topological invariants and applications were discussed in
~\cite{E09}. The articles ~\cite{GS05} and ~\cite{RS01} are crucial for the
understanding of the Fredholm and properness properties of the quasilinear
elliptic systems of the second order and of the operators of this kind
on ${\mathbb R}^{N}$.
The non-Fredholm operators appear also when considering the wave systems with
an infinite number of localized traveling waves (see ~\cite{AMP14}). In
particular, when $a$ vanishes, the operator $A$ is Fredholm in certain properly chosen
weighted spaces (see \cite{Amrouche1997}, \cite{Amrouche2008},
\cite{Bolley1993}, \cite{Bolley2001}, \cite{B88}). However, the situation when
$a \neq 0$ is considerably
different and the methods developed in these works cannot be applied.
The front propagation problems in the context of the anomalous diffusion were studied actively in
~\cite{VNN10}, ~\cite{VNN13}.

\medskip

We set $K_{m}(x) = \varepsilon_{m} H_{m}(x)$, with
$\varepsilon_{m} \geq 0$, such that
\begin{equation}
\label{es}  
\varepsilon:=\hbox{max}_{1\leq m\leq N}\varepsilon_{m}, \quad s_{1}:=\hbox{min}
_{1\leq m\leq N} s_{1, m}, \quad S_{1}:=\hbox{max}_{1\leq m\leq N} s_{1, m},
\end{equation}
where $\displaystyle{\frac{1}{4}<s_{1}, \ S_{1}<\frac{3}{4}}$. Let us assume
the following.

\bigskip

\noindent
{\bf Assumption 1.1.}  {\it Let
$\displaystyle{1\leq m\leq N, \ \frac{1}{4}<s_{1, m}<\frac{3}{4}, \
s_{1, m}<s_{2, m}<1}$, the functions   
$f_{m}(x): {\mathbb R}^{3}\to {\mathbb R}$ are nontrivial for a
certain $m$, so that
$$
f_{m}(x)\in L^{1}({\mathbb R}^{3}), \quad
(-\Delta)^{1-s_{1, m}} f_{m}(x)\in L^{2}({\mathbb R}^{3}).
$$
Additionally, we assume that
$H_{m}(x): {\mathbb R}^{3}\to {\mathbb R}$, such that
$$
H_{m}(x)\in L^{1}({\mathbb R}^{3}), \quad
(-\Delta)^{1-s_{1, m}}H_{m}(x)\in L^{2}({\mathbb R}^{3}).
$$
Furthermore,
\begin{equation}
\label{h2}  
H^{2}:=\sum_{m=1}^{N}\|H_{m}(x)\|_{L^{1}({\mathbb R}^{3})}^{2}>0
\end{equation}
and}
\begin{equation}
\label{q2}  
Q^{2}:=\sum_{m=1}^{N}\|(-\Delta)^{1-s_{1, m}}
H_{m}(x)\|_{L^{2}({\mathbb R}^{3})}^{2}>0.
\end{equation}

\bigskip

We work here in the space of dimension $d=3$. This is related to the
solvability relations for the linear Poisson type problem (\ref{lp}) formulated in
Lemma 4.1 further down. From the point of view of the practical applications, the space dimension is not
limited to $d=3$, since the space variable here is correspondent to the cell
genotype but not to the usual physical space. 
We apply the Sobolev inequality for the fractional negative Laplacian
(see Lemma 2.2 of ~\cite{HYZ12}, also ~\cite{L83}), namely
\begin{equation}
\label{sobfr}
\|f_{m}\|_{L^{\frac{6}{4s_{1, m}-1}}({\mathbb R}^{3})}\leq c_{s_{1, m}}
\|(-\Delta)^{1-s_{1, m}}f_{m}\|_{L^{2}({\mathbb R}^{3})},    
\end{equation}
where $\displaystyle{\frac{1}{4}<1-s_{1, m}<\frac{3}{4}}$ and
$1\leq m\leq N$. 
By means of  Assumption 1.1 above and the standard interpolation
argument, we obtain
\begin{equation}
\label{l2}
f_{m}\in L^{2}({\mathbb R}^{3}), \quad 1\leq m\leq N.
\end{equation}
For the technical purposes, we use the Sobolev spaces 
\begin{equation}
\label{h2sm}  
H^{2s_{2, m}}({\mathbb R}^{3}):=\{\phi:{\mathbb R}^{3}\to {\mathbb {\mathbb R}}
\ | \ \phi\in L^{2}({\mathbb R}^{3}), \ (-\Delta)^{s_{2, m}} \phi \in
L^{2}({\mathbb R}^{3}) \}, 
\end{equation}
with
$\displaystyle{0<s_{2, m}\leq 1, \ 1\leq m\leq N}$. Each space (\ref{h2sm}) is equipped with the norm
\begin{equation}
\label{n}
\|\phi\|_{H^{2s_{2, m}}({\mathbb R}^{3})}^{2}:=\|\phi\|_{L^{2}({\mathbb R}^{3})}^{2}+
\|(-\Delta)^{s_{2, m}} \phi \|_{L^{2}({\mathbb R}^{3})}^{2}.
\end{equation}
For a vector function
$$
u(x)=(u_{1}(x), u_{2}(x), ..., u_{N}(x))^{T},
$$
in our article we will use the norms
\begin{equation}
\label{u1N}
\|u\|_{H^{2}({\mathbb R}^{3}, {\mathbb R}^{N})}^{2}:=\sum_{m=1}^{N}\|u_{m}\|_{H^{2}({\mathbb R}^{3}) }^{2}=
\sum_{m=1}^{N}\{\|u_{m}\|_{L^{2}({\mathbb R}^{3}) }^{2}+\|\Delta u_{m}\|_{L^{2}({\mathbb R}^{3})}^{2}\}, 
\end{equation}
and
$$
\|u\|_{L^{2}({\mathbb R}^{3}, {\mathbb R}^{N})}^{2}:=\sum_{m=1}^{N}
\|u_{m}\|_{L^{2}({\mathbb R}^{3})}^{2}.
$$
According to the standard Sobolev embedding,
\begin{equation}
\label{e}
\|\phi\|_{L^{\infty}({\mathbb R}^{3})}\leq c_{e}\|\phi\|_{H^{2}({\mathbb R}^{3})},
\end{equation}
where $c_{e}>0$ is the constant of the embedding.
When all the nonnegative parameters $\varepsilon_{m}$ vanish, we arrive at
the linear Poisson type equations
\begin{equation}
\label{lpm}
[(-\Delta)^{s_{1, m}}+(-\Delta)^{s_{2, m}}]u_{m}(x)=f_{m}(x), \quad 1\leq m\leq N.
\end{equation}
By means of Lemma 4.1 below under the given conditions each problem
(\ref{lpm}) has a unique solution
$$
u_{0, m}\in H^{2s_{2, m}}({\mathbb R}^{3}), \quad
\frac{1}{4}<s_{1, m}<\frac{3}{4}, \quad s_{1, m}<s_{2, m}<1, \quad 1\leq m\leq N,
$$
and no orthogonality relations for the right side of (\ref{lpm}) are required
here. Evidently, for $1\leq m\leq N$,
\begin{equation}
\label{ds2ms1m} 
[-\Delta+(-\Delta)^{1+s_{2, m}-s_{1, m}}]u_{0,m}(x)=
(-\Delta)^{1-s_{1, m}}f_{m}(x)\in
L^{2}({\mathbb R}^{3})
\end{equation}
due to Assumption 1.1. It can be trivially obtained from (\ref{ds2ms1m}) via the
standard Fourier transform (\ref{f}) that
$$
\Delta u_{0,m}\in L^{2}({\mathbb R}^{3}), \ 1\leq m\leq N.
$$
Thus, each linear problem (\ref{lpm}) admits a unique solution
$u_{0,m}\in H^{2}({\mathbb R}^{3})$. According to the definition of the norm
(\ref{u1N}), 
$$
u_{0}(x):=(u_{0,1}(x), u_{0,2}(x), ..., u_{0,N}(x))^{T}\in
H^{2}({\mathbb R}^{3}, {\mathbb R}^{N}).
$$
We seek the resulting solution of the nonlinear system of equations
(\ref{p}) as
\begin{equation}
\label{r}
u(x)=u_{0}(x)+u_{p}(x)
\end{equation}
with
$$
u_{p}(x):=(u_{p,1}(x), u_{p,2}(x),...,u_{p,N}(x))^{T}.
$$
Obviously, we arrive at the perturbative system of equations
\begin{equation}
\label{pert}
[(-\Delta)^{s_{1, m}}+(-\Delta)^{s_{2, m}}]u_{p,m}(x)=\varepsilon_{m} \int_{{\mathbb R}^{3}}
H_{m}(x-y)g_{m}(u_{0}(y)+u_{p}(y))dy,
\end{equation}
where
$\displaystyle{\frac{1}{4}<s_{1, m}<\frac{3}{4}, \ s_{1, m}<s_{2, m}<1, \
1\leq m\leq N}$.

Let us  introduce a closed ball in the Sobolev space
\begin{equation}
\label{b}
B_{\rho}:=\{u\in H^{2}({\mathbb R}^{3}, {\mathbb R}^{N}) \ | \ \|u\|_
{H^{2}({\mathbb R}^{3}, {\mathbb R}^{N})}\leq \rho \}, \quad 0<\rho\leq 1.
\end{equation}
We look for the solution of system (\ref{pert}) as the fixed point of the
auxiliary nonlinear problem
\begin{equation}
\label{aux}
[(-\Delta)^{s_{1, m}}+(-\Delta)^{s_{2, m}}]u_{m}(x)=\varepsilon_{m} \int_{{\mathbb R}^{3}}
H_{m}(x-y)g_{m}(u_{0}(y)+v(y))dy,
\end{equation}
with
$\displaystyle{\frac{1}{4}<s_{1, m}<\frac{3}{4}, \ s_{1, m}<s_{2, m}<1, \
1\leq m\leq N}$
in ball (\ref{b}).

 For a
given vector function $v(y)$ this is a system of equations with respect to
$u(x)$.
The left side of the $m$th equation in (\ref{aux}) contains the operator which
does not satisfy the Fredholm property
\begin{equation}
\label{sm}  
l_{m}:=(-\Delta)^{s_{1, m}}+(-\Delta)^{s_{2, m}}: H^{2s_{2, m}}({\mathbb R}^{3})\to
L^{2}({\mathbb R}^{3}), \quad 1\leq m\leq N.
\end{equation}
(\ref{sm}) is defined by means of the spectral calculus. It is the
pseudo-differential operator with the symbol $|p|^{2s_{1, m}}+|p|^{2s_{2, m}}$,
so that for $1\leq m\leq N$
$$
l_{m}\phi(x)=\frac{1}{(2\pi)^{\frac{3}{2}}}\int_{{\mathbb R}^{3}}
(|p|^{2s_{1, m}}+|p|^{2s_{2, m}})\widehat{\phi}(p)e^{ipx}dp, \quad
\phi\in H^{2s_{2, m}}({\mathbb R}^{3}), 
$$
where the standard Fourier transform is introduced in (\ref{f}).
The essential spectrum of (\ref{sm}) fills the nonnegative semi-axis
$[0, +\infty)$.
Hence, such operator does not have a bounded inverse. The analogous situation
appeared in works ~\cite{VV111} and ~\cite{VV14} but as distinct from the
present case, the problems treated there required the orthogonality
relations.  Persistence of pulses for certain local reaction-diffusion equations via
the fixed point technique was covered in ~\cite{CV21}.
But the Schr\"odinger type operator contained in the nonlinear
problem there possessed the Fredholm property.

We introduce the closed ball in the
space of $N$ dimensions as
\begin{equation}
\label{i}
I:=\{z\in {\mathbb R}^{N} \ | \ |z|_{{\mathbb R}^{N}}\leq
c_{e}\|u_{0}\|_{H^{2}({\mathbb R}^{3}, {\mathbb R}^{N})}+c_{e} \}.
\end{equation}
Here and further down $|.|_{{\mathbb R}^{N}}$ will stand for the length of a vector in
${\mathbb R}^{N}$.
The closed ball $D_{M}$ in the space of $C^{2}(I,{\mathbb R}^{N})$
vector functions is 
\begin{equation}
\label{M}
\{g(z):=(g_{1}(z), g_{2}(z),..., g_{N}(z))\in C^{2}(I,{\mathbb R}^{N}) \ | \
\|g\|_{C^{2}(I,{\mathbb R}^{N})}\leq M \}
\end{equation}
with  $M>0$.
In this context the norms
\begin{equation}
\label{gN}
\|g\|_{C^{2}(I,{\mathbb R}^{N})}:=\sum_{m=1}^{N}\|g_{m}\|_{C^{2}(I)},
\end{equation}
\begin{equation}
\label{gn}
\|g_{m}\|_{C^{2}(I)}:=\|g_{m}\|_{C(I)}+\sum_{n=1}^{N}\Big\|\frac{\partial g_{m}}
{\partial z_{n}}\Big\|_{C(I)}+\sum_{n,l=1}^{N}\Big\|\frac{\partial^{2}g_{m}}
{\partial z_{n}\partial z_{l}}\Big\|_{C(I)},
\end{equation}
where $\|g_{m}\|_{C(I)}:=\hbox{max}_{z\in I}|g_{m}(z)|$.
Let us impose the following technical conditions on the nonlinear part of system
(\ref{p}). From the point of view of the applications in biology,
$g_{m}(z)$ can be, for instance the quadratic functions,  describing the
cell-cell interactions.

\bigskip

\noindent
{\bf Assumption 1.2.} {\it Let $1\leq m\leq N$. We assume that
$g_{m}: {\mathbb R}^{N}\to {\mathbb R}$ is such that
$g_{m}(0)=0$ and $\nabla g_{m}(0)=0$. Let us also suppose that $g\in D_{M}$ and
it does not vanish identically in the ball $I$}.

\bigskip

Let us use the auxiliary Assumptions 1.1 and 1.2 above in the proofs of our main
statements. It is not evident at the moment if there is a more efficient way to
analyse our system  which would allow us to weaken such
conditions.

We introduce the operator $\tau_g$, so that $u = \tau_g v$, where $u$ is a
solution of the system of equations (\ref{aux}). Our first main result is
as follows.

\bigskip

\noindent
{\bf Theorem 1.3.} {\it Let Assumptions 1.1 and 1.2 be valid. Then for every
$\rho\in (0, 1]$ problem
(\ref{aux}) defines the map $\tau_{g}: B_{\rho}\to B_{\rho}$, which is a strict
contraction for all
$$
0<\varepsilon\leq
\frac{\rho}{M(\|u_{0}\|_{H^{2}({\mathbb R}^{3}, {\mathbb R}^{N})}+1)^{2}}\times
$$
\begin{equation}
\label{eps}
\Bigg[H^{2}(\|u_{0}\|_{H^{2}({\mathbb R}^{3}, {\mathbb R}^{N})}+1)^{\frac{8S_{1}}{3}-2}\frac{1}{(4s_{1})^{\frac{4s_{1}}{3}}}
\frac{3}{(3-4S_{1})(2\pi^{2})^{\frac{4s_{1}}{3}}}   
+Q^{2}
\Bigg]^{-\frac{1}{2}}.
\end{equation}
The unique fixed point $u_{p}$ of
such map $\tau_{g}$ is the only solution of system (\ref{pert}) in $B_{\rho}$.}

\bigskip

The constants $\varepsilon, \ s_{1}, \ S_{1}\ H, \ Q$ in formula (\ref{eps}) are introduced in 
(\ref{es}),  (\ref{h2}) and (\ref{q2}).

Note that the resulting solution $u(x)$ of the system of equations (\ref{p})
given by (\ref{r}) will be nontrivial since
the influx/efflux terms $f_{m}(x)$ do not vanish identically for a certain
$1\leq m\leq N$ and all $g_{m}(0)=0$ as assumed. We have
the following trivial lemma.

\bigskip

\noindent
{\bf Lemma 1.4.} {\it Let $R\in (0, +\infty)$. Consider the
function
$$
\varphi(R):=\alpha R^{3-4s}+\frac{1}{R^{4s}}, \quad
\frac{1}{4}<s<\frac{3}{4}, \quad \alpha>0.
$$
It achieves its minimal value at \
$\displaystyle{R^{*}:=\Bigg({\frac{4s}{\alpha (3-4s)}}\Bigg)^{\frac{1}{3}}}$,
which is given by}
$$
\varphi(R^{*})=3(3-4s)^{\frac{4s}{3}-1}(4s)^{-\frac{4s}{3}}\alpha^{\frac{4s}{3}}.
$$

\bigskip

The second main statement of the article is devoted to the continuity of the resulting
solution of system  (\ref{p}) given by formula (\ref{r}) with
respect to the nonlinear vector function $g$. We define the following
positive technical expression
$$
\sigma:=M(\|u_{0}\|_{H^{2}({\mathbb R}^{3}, {\mathbb R}^{N})}+1)\times
$$
\begin{equation}
\label{sig}
\Bigg\{\frac{H^{2}(\|u_{0}\|
_{H^{2}({\mathbb R}^{3}, {\mathbb R}^{N})}+1)^{\frac{8S_{1}}{3}-2}3}
{(3-4S_{1}){(2{\pi}^{2})}^{\frac{4s_{1}}{3}}}\frac{1}{{(4s_{1})}^{\frac{4s_{1}}{3}}}+
Q^{2}\Bigg\}^{\frac{1}{2}}.
\end{equation}

\bigskip

\noindent
{\bf Theorem 1.5.} {\it Let $j=1,2$, the conditions of Theorem 1.3 hold,
so that
$u_{p,j}$ is the unique fixed point of the map
$\tau_{g_{j}}: B_{\rho}\to B_{\rho}$, which is a strict contraction for all the values
of $\varepsilon$ satisfying inequality (\ref{eps}) and the resulting solution of the
system of equations (\ref{p}) with $g(z)=g_{j}(z)$ equals to
\begin{equation}
\label{jres}
u_{j}(x):=u_{0}(x)+u_{p, j}(x).
\end{equation}
Then for all the values of $\varepsilon$, which satisfy bound (\ref{eps}),
the estimate
\begin{equation}
\label{cont}
\|u_{1}-u_{2}\|_{H^{2}({\mathbb R}^{3}, {\mathbb R}^{N})}\leq \frac{\varepsilon \sigma}
{M(1-\varepsilon \sigma)}(\|u_{0}\|_{H^{2}({\mathbb R}^{3}, {\mathbb R}^{N})}+1)
\|g_{1}-g_{2}\|_{C^{2}(I, {\mathbb R}^{N})}
\end{equation}
is valid.}

\bigskip

Let us proceed to proving our first main proposition.

\bigskip

%%%%%%%%%%%%%%%%%%%%%%%%%%%%%%%%%%%%%%%%%%%%%%%%%%

\setcounter{section}{2}
\setcounter{equation}{0}

\centerline{\bf 2. The existence of the perturbed solution}

\bigskip

\noindent
{\it Proof of Theorem 1.3.} We choose an arbitrary  vector function
$v\in B_{\rho}$ and denote the terms involved in the integral expressions
in the right side of system  (\ref{aux}) as
$$
G_{m}(x):=g_{m}(u_{0}(x)+v(x)), \quad 1\leq m\leq N.
$$
The standard Fourier transform is given by
\begin{equation}
\label{f}
\widehat{\phi}(p):=\frac{1}{(2\pi)^{\frac{3}{2}}}\int_{{\mathbb R}^{3}}\phi(x)e^{-ipx}dx,
\quad p\in {\mathbb R}^{3}.
\end{equation}
Evidently, the upper bound
\begin{equation}
\label{fub}
\|\widehat{\phi}\|_{L^{\infty}({\mathbb R}^{3})}\leq \frac{1}{(2\pi)^{\frac{3}{2}}}
\|\phi\|_{L^{1}({\mathbb R}^{3})}
\end{equation}
holds.
Let us apply (\ref{f}) to both sides of the system of equations (\ref{aux}). This yields
$$
\widehat{u_{m}}(p)=\varepsilon_{m} (2\pi)^{\frac{3}{2}}
\frac{\widehat{H_{m}}(p)\widehat{G_{m}}(p)}{|p|^{2s_{1, m}}+|p|^{2s_{2, m}}}
$$
with
$\displaystyle{\frac{1}{4}<s_{1, m}<\frac{3}{4}, \ s_{1, m}<s_{2, m}<1, \
1\leq m\leq N}$.
Thus, we have the expression for the norm as
\begin{equation}
\label{un}
\|u_{m}\|_{L^{2}({\mathbb R}^{3})}^{2}={(2\pi)}^{3}{\varepsilon}_{m}^{2}\int_{{\mathbb R}^{3}}
\frac{|\widehat{H_{m}}(p)|^{2}|\widehat{G_{m}}(p)|^{2}}
{[|p|^{2s_{1, m}}+|p|^{2s_{2, m}}]^{2}}dp.
\end{equation}
As distinct from the previous articles ~\cite{VV111} and ~\cite{VV14} containing the standard
Laplace operator in the diffusion term, here we do not try to control the norms
$$
\Bigg\|\frac{\widehat{H_{m}}(p)}{|p|^{2s_{1, m}}+|p|^{2s_{2, m}}}\Bigg\|_
{L^{\infty}({\mathbb R}^{3})}, \quad 1\leq m\leq N.
$$
Instead, we estimate the right side of (\ref{un}) via the analog of
inequality (\ref{fub}) applied to functions $H_{m}$ and $G_{m}$ with
$R\in (0, +\infty)$ as
$$
{(2\pi)}^{3}{\varepsilon}_{m}^{2}\int_{{\mathbb R}^{3}}
\frac{|\widehat{H_{m}}(p)|^{2}|\widehat{G_{m}}(p)|^{2}}
{[|p|^{2s_{1, m}}+|p|^{2s_{2, m}}]^{2}}dp\leq
$$
$$
{(2\pi)}^{3}{\varepsilon}_{m}^{2}\Big[\int_{|p|\leq R}\frac{|\widehat{H_{m}}(p)|
^{2}|\widehat{G_{m}}(p)|^{2}}{|p|^{4s_{1, m}}}dp+\int_{|p|>R}\frac{|\widehat{H_{m}}
(p)|^{2}|\widehat{G_{m}}(p)|^{2}}{|p|^{4s_{1, m}}}dp \Big]\leq
$$
\begin{equation}
\label{ub1}
 {\varepsilon}_{m}^{2}\|H_{m}\|_{L^{1}({\mathbb R}^{3})}^{2}\Bigg\{
\frac{1}{2{\pi}^{2}}\|G_{m}\|_{L^{1}({\mathbb R}^{3})}^{2}\frac{R^{3-4s_{1, m}}}
{3-4s_{1, m}}+
\frac{\|G_{m}\|_{L^{2}({\mathbb R}^{3})}^{2}}{R^{4s_{1, m}}}\Bigg\}.
\end{equation}
Let us recall the norm definition (\ref{u1N}). Using the triangle
inequality along with the fact that $v\in B_{\rho}$, we easily obtain
$$
\|u_{0}+v\|_{L^{2}({\mathbb R}^{3}, {\mathbb R}^{N})}\leq
\|u_{0}\|_{H^{2}({\mathbb R}^{3},{\mathbb R}^{N})}+1.
$$
Sobolev embedding (\ref{e}) implies that
$$
|u_{0}+v|_{{\mathbb R}^{N}}\leq c_{e}(\|u_{0}\|_{H^{2}({\mathbb R}^{3}, {\mathbb R}^{N})}+1).
$$
Obviously,
$$
G_{m}(x)=\int_{0}^{1}\nabla g_{m}(t(u_{0}(x)+v(x))).(u_{0}(x)+v(x))dt, \quad
1\leq m\leq N.
$$
Here and below the dot will denote the scalar product of two vectors in ${\mathbb R}^{N}$.
We use the ball $I$ introduced in (\ref{i}). Hence,
$$
|G_{m}(x)|\leq \hbox{sup}_{z\in I}|\nabla g_{m}(z)|_{{\mathbb R}^{N}}
|u_{0}(x)+v(x)|_{{\mathbb R}^{N}}\leq M|u_{0}(x)+v(x)|_{{\mathbb R}^{N}}.
$$
Then
$$
\|G_{m}\|_{L^{2}({\mathbb R}^{3})}\leq M\|u_{0}+v\|_{L^{2}({\mathbb R}^{3},{\mathbb R}^{N})}
\leq M(\|u_{0}\|_{H^{2}({\mathbb R}^{3},{\mathbb R}^{N})}+1).
$$
Note that for $t\in [0,1]$ and $1\leq m,j\leq N$,
$$
\frac{\partial g_{m}}{\partial z_{j}}(t(u_{0}(x)+v(x)))=\int_{0}^{t}\nabla
\frac{\partial g_{m}}{\partial z_{j}}(\tau(u_{0}(x)+v(x))).(u_{0}(x)+v(x))d\tau.
$$
Therefore,
$$
\Big|\frac{\partial g_{m}}{\partial z_{j}}(t(u_{0}(x)+v(x)))\Big|\leq
\hbox{sup}_{z\in I}\Big|\nabla \frac{\partial g_{m}}{\partial z_{j}}\Big|_
{{\mathbb R}^{N}}|u_{0}(x)+v(x)|_{{\mathbb R}^{N}}\leq
$$
$$
\sum_{n=1}^{N}\Big\|\frac{\partial^{2} g_{m}}{\partial z_{n}\partial z_{j}}
\Big\|_{C(I)}|u_{0}(x)+v(x)|_{{\mathbb R}^{N}}.
$$
This means that
$$
|G_{m}(x)|\leq |u_{0}(x)+v(x)|_{{\mathbb R}^{N}}\sum_{j=1}^{N}
\sum_{n=1}^{N}\Big\|\frac{\partial^{2} g_{m}}
{\partial z_{n}\partial z_{j}}\Big\|_{C(I)}|u_{0,j}(x)+v_{j}(x)|\leq
$$
$$
M|u_{0}(x)+v(x)|_{{\mathbb R}^{N}}^{2}.
$$
Hence,
\begin{equation}
\label{G1}
\|G_{m}\|_{L^{1}({\mathbb R}^{3})}\leq M\|u_{0}+v\|_{L^{2}({\mathbb R}^{3},
{\mathbb R}^{N})}^{2}\leq M(\|u_{0}\|_{H^{2}({\mathbb R}^{3}, {\mathbb R}^{N})}+1)^{2}.
\end{equation}
This enables us to obtain the estimate from above for the right side of
(\ref{ub1}) equal to
$$
{\varepsilon}_{m}^{2}M^{2}\|H_{m}\|_{L^{1}({\mathbb R}^{3})}^{2}
(\|u_{0}\|_{H^{2}({\mathbb R}^{3}, {\mathbb R}^{N})}+1)^{2}\times
$$
$$
\Bigg\{\frac{(\|u_{0}\|_{H^{2}({\mathbb R}^{3}, {\mathbb R}^{N})}+
1)^{2}R^{3-4s_{1, m}}}{2{\pi}^{2}(3-4s_{1, m})}
+\frac{1}{R^{4s_{1, m}}}\Bigg\},
$$
where $R\in (0, +\infty)$. Lemma 1.4 gives us the minimal value
of the quantity above, so that
$$
\|u_{m}\|_{L^{2}({\mathbb R}^{3})}^{2}\leq {\varepsilon_{m}}^{2}M^{2}\|H_{m}\|_
{L^{1}({\mathbb R}^{3})}^{2}\times
$$
$$
(\|u_{0}\|_{H^{2}({\mathbb R}^{3}, {\mathbb R}^{N})}+1)^
{2+\frac{8s_{1, m}}{3}}
\frac{1}{(4s_{1, m})^{\frac{4s_{1, m}}{3}}}\frac{3}
{(3-4s_{1, m})(2{\pi}^{2})^{\frac{4s_{1, m}}{3}}}.
$$
Let us recall definition (\ref{es}). Thus, for $1\leq m\leq N$
$$
\|u_{m}\|_{L^{2}({\mathbb R}^{3})}^{2}\leq {\varepsilon}^{2}M^{2} \|H_{m}\|_{L^{1}({\mathbb R}^{3})}^{2}\times
$$
$$
(\|u_{0}\|_{H^{2}({\mathbb R}^{3}, {\mathbb R}^{N})}+1)^
{2+\frac{8S_{1}}{3}}
\frac{1}{(4s_{1})^{\frac{4s_{1}}{3}}}\frac{3}
{(3-4S_{1})(2{\pi}^{2})^{\frac{4s_{1}}{3}}},
$$
such that
$$
\|u\|_{L^{2}({\mathbb R}^{3},{\mathbb R}^{N})}^{2}\leq {\varepsilon}^{2}M^{2}H^{2}
\times
$$
\begin{equation}
\label{ul2ub}
(\|u_{0}\|_{H^{2}({\mathbb R}^{3}, {\mathbb R}^{N})}+1)^
{2+\frac{8S_{1}}{3}}
\frac{1}{(4s_{1})^{\frac{4s_{1}}{3}}}\frac{3}
{(3-4S_{1})(2{\pi}^{2})^{\frac{4s_{1}}{3}}}.
\end{equation}
By virtue of (\ref{aux}),
$$
[-\Delta+(-\Delta)^{1+s_{2, m}-s_{1, m}}]u_{m}(x)=\varepsilon_{m}
(-\Delta)^{1-s_{1, m}}
\int_{{\mathbb R}^{3}}H_{m}(x-y)G_{m}(y)dy,
$$
where
$\displaystyle{\frac{1}{4}<s_{1,m}<\frac{3}{4}, \ s_{1, m}<s_{2, m}<1, \
1\leq m\leq N}$.

Let us apply here the standard Fourier transform (\ref{f}),
the analog of estimate (\ref{fub}) used for function $G_{m}$
and (\ref{G1}). This yields
$$
\|\Delta u_{m}\|_{L^{2}({\mathbb R}^{3})}^{2}\leq \varepsilon_{m}^{2}
\|G_{m}\|_{L^{1}({\mathbb R}^{3})}^{2}
\|(-\Delta)^{1-s_{1, m}}H_{m}\|_
{L^{2}({\mathbb R}^{3})}^{2}\leq
$$
$$
\varepsilon^{2}M^{2}
(\|u_{0}\|_{H^{2}({\mathbb R}^{3}, {\mathbb R}^{N})}+1)^{4}
\|(-\Delta)^{1-s_{1, m}}H_{m}\|_{L^{2}({\mathbb R}^{3})}^{2}.
$$
Hence,
\begin{equation}
\label{32}
\sum_{m=1}^{N}\|\Delta u_{m}\|_{L^{2}({\mathbb R}^{3})}^{2}\leq
\varepsilon^{2}M^{2}(\|u_{0}\|_{H^{2}({\mathbb R}^{3}, {\mathbb R}^{N})}+1)^{4}Q^{2}.
\end{equation}
We recall the definition of the norm (\ref{u1N}). Inequalities
(\ref{ul2ub}) and (\ref{32}) imply that
$$
\|u\|_{H^{2}({\mathbb R}^{3}, {\mathbb R}^{N})}\leq
\varepsilon M(\|u_{0}\|_{H^{2}({\mathbb R}^{3}, {\mathbb R}^{N})}+1)^{2}\times
$$
\begin{equation}
\label{rh}
\Bigg[H^{2}(\|u_{0}\|_{H^{2}({\mathbb R}^{3}, {\mathbb R}^{N})}+1)^{\frac{8S_{1}}{3}-2}
\frac{1}{(4s_{1})^{\frac{4s_{1}}{3}}}\frac{3}
{(3-4S_{1})(2{\pi}^{2})^{\frac{4s_{1}}{3}}}
+Q^{2}\Bigg]
^{\frac{1}{2}}
\leq \rho
\end{equation}
for all the values of $\varepsilon$ satisfying condition (\ref{eps}). Therefore,
$u\in B_{\rho}$ as well.

Let us suppose that for some $v\in B_{\rho}$ the system of equations (\ref{aux}) admits
two solutions $u_{1,2}\in B_{\rho}$. Evidently, their
difference $w(x):=u_{1}(x)-u_{2}(x) \in H^{2}({\mathbb R}^{3}, {\mathbb R}^{N})$
solves the homogeneous system 
$$
[(-\Delta)^{s_{1, m}}+(-\Delta)^{s_{2, m}}]w_{m}(x)=0
$$
with
$\displaystyle{\frac{1}{4}<s_{1, m}<\frac{3}{4}, \ s_{1, m}<s_{2, m}<1, \
1\leq m\leq N}$.

Note that each operator $l_{m}: H^{2s_{2, m}}({\mathbb R}^{3})\to L^{2}({\mathbb R}^{3})$
defined in (\ref{sm}) does not possess any nontrivial zero modes. Hence,
$w(x)\equiv 0$ in ${\mathbb R}^{3}$. This means that problem
(\ref{aux}) defines a map
$\tau_{g}: B_{\rho}\to B_{\rho}$ for all $\varepsilon$, which satisfy assumption
(\ref{eps}).

The goal is to demonstrate that such map is a strict contraction. We choose
arbitrarily $v_{1}, v_{2}\in B_{\rho}$. By means of the reasoning above, 
$u_{1,2}:=\tau_{g}v_{1,2}\in B_{\rho}$ as well for $\varepsilon$ satisfying (\ref{eps}).
Clearly, according to (\ref{aux}) we have
\begin{equation}
\label{aux1}
[(-\Delta)^{s_{1, m}}+(-\Delta)^{s_{2, m}}]u_{1,m}(x)=\varepsilon_{m} \int_{{\mathbb R}^{3}}
H_{m}(x-y)g_{m}(u_{0}(y)+v_{1}(y))dy,
\end{equation}
\begin{equation}
\label{aux2}
[(-\Delta)^{s_{1, m}}+(-\Delta)^{s_{2, m}}]u_{2,m}(x)=\varepsilon_{m} \int_{{\mathbb R}^{3}}
H_{m}(x-y)g_{m}(u_{0}(y)+v_{2}(y))dy,
\end{equation}
where $\displaystyle{\frac{1}{4}<s_{1, m}<\frac{3}{4}, \ s_{1, m}<s_{2, m}<1, \ 1\leq m\leq N}$.
 We introduce
$$
G_{1, m}(x):=g_{m}(u_{0}(x)+v_{1}(x)), \quad G_{2, m}(x):=g_{m}(u_{0}(x)+v_{2}(x)),
\quad 1\leq m\leq N.
$$
Let us apply the standard Fourier transform (\ref{f}) to both sides of
systems of equations (\ref{aux1}) and (\ref{aux2}). This yields
$$
\widehat{u_{1, m}}(p)=\varepsilon_{m} (2\pi)^{\frac{3}{2}}
\frac{\widehat{H_{m}}(p)\widehat{G_{1, m}}(p)}{|p|^{2s_{1, m}}+|p|^{2s_{2, m}}}, \quad
\widehat{u_{2, m}}(p)=\varepsilon_{m} (2\pi)^{\frac{3}{2}}
\frac{\widehat{H_{m}}(p)\widehat{G_{2, m}}(p)}{|p|^{2s_{1, m}}+|p|^{2s_{2, m}}}.
$$
Obviously,
\begin{equation}
\label{u12mn}  
\|u_{1, m}-u_{2, m}\|_{L^{2}({\mathbb R}^{3})}^{2}=\varepsilon_{m}^{2}{(2\pi)}^{3}
\int_{{\mathbb R}^{3}} \frac{|\widehat{H_{m}}(p)|^{2}
|{\widehat{G_{1, m}}(p)}-{\widehat{G_{2, m}}}(p)|^{2}}
{[|p|^{2s_{1, m}}+|p|^{2s_{2, m}}]^{2}}dp.
\end{equation}
Evidently, the right side of (\ref{u12mn}) can be bounded from above
by virtue of estimate (\ref{fub}) as
$$
\varepsilon_{m}^{2}{(2\pi)}^{3}\Bigg[\int_{|p|\leq R}\frac{|\widehat{H_{m}}(p)|^{2}
|{\widehat{G_{1, m}}(p)}-{\widehat{G_{2, m}}}(p)|^{2}}{|p|^{4s_{1, m}}}dp+
$$
$$  
\int_{|p|>R}\frac{|\widehat{H_{m}}(p)|^{2}
|{\widehat{G_{1, m}}(p)}-{\widehat{G_{2, m}}}(p)|^{2}}{|p|^{4s_{1, m}}}dp \Bigg]\leq  
\varepsilon^{2} \|H_{m}\|_{L^{1}({\mathbb R}^{3})}^{2} \times
$$
$$
\Bigg\{\frac {\|G_{1, m}-G_{2, m}\|_{L^{1}({\mathbb R}^{3})}^{2}}
{2{\pi}^{2}}\frac{R^{3-4s_{1, m}}}{3-4s_{1, m}}+
\frac{\|G_{1, m}-G_{2, m}\|_{L^{2}({\mathbb R}^{3})}^{2}}{R^{4s_{1, m}}}\Bigg\},
$$
where $R\in (0,+\infty)$. Note that for
$1\leq m\leq N$, we have
$$
G_{1, m}(x)-G_{2, m}(x)=\int_{0}^{1}\nabla g_{m}(u_{0}(x)+tv_{1}(x)+(1-t)v_{2}(x)).
(v_{1}(x)-v_{2}(x))dt.
$$
Clearly, for $t\in [0,1]$
$$
\|v_{2}+t(v_{1}-v_{2})\|_{H^{2}({\mathbb R}^{3}, {\mathbb R}^{N})}\leq
t\|v_{1}\|_{H^{2}({\mathbb R}^{3}, {\mathbb R}^{N})}+
$$
$$
(1-t)\|v_{2}\|_{H^{2}({\mathbb R}^{3}, {\mathbb R}^{N})}
\leq \rho.
$$
Therefore, $v_{2}+t(v_{1}-v_{2})\in B_{\rho}$. Let us derive the upper
bound
$$
|G_{1, m}(x)-G_{2, m}(x)|\leq \hbox{sup}_{z\in I}|\nabla g_{m}(z)|_{{\mathbb R}^{N}}
|v_{1}(x)-v_{2}(x)|_{{\mathbb R}^{N}}
\leq M|v_{1}(x)-v_{2}(x)|_{{\mathbb R}^{N}},
$$
such that
$$
\|G_{1, m}-G_{2, m}\|_{L^{2}({\mathbb R}^{3})}\leq M\|v_{1}-v_{2}\|_
{L^{2}({\mathbb R}^{3}, {\mathbb R}^{N})}\leq M\|v_{1}-v_{2}\|_
{H^{2}({\mathbb R}^{3}, {\mathbb R}^{N})}.
$$
We express
$\displaystyle{\frac{\partial g_{m}}{\partial z_{j}}(u_{0}(x)+tv_{1}(x)+(1-t)
v_{2}(x))}$  for $1\leq m,j\leq N$ as
$$
\int_{0}^{1}\nabla \frac{\partial g_{m}}{\partial z_{j}}
(\tau[u_{0}(x)+tv_{1}(x)+(1-t)v_{2}(x)]).[u_{0}(x)+tv_{1}(x)+(1-t)v_{2}(x)]d\tau.
$$
Hence, for $t\in [0,1]$
$$
\Big|\frac{\partial g_{m}}{\partial z_{j}}(u_{0}(x)+tv_{1}(x)+(1-t)v_{2}(x))\Big|\
\leq
$$
$$
\leq\sum_{n=1}^{N}\Bigg\|\frac{\partial^{2}g_{m}}{\partial z_{n}\partial z_{j}}
\Bigg\|_{C(I)}(|u_{0}(x)|_{{\mathbb R}^{N}}+t|v_{1}(x)|_{{\mathbb R}^{N}}+
(1-t)|v_{2}(x)|_{{\mathbb R}^{N}}).
$$
Then
$$
|G_{1, m}(x)-G_{2, m}(x)|\leq 
M|v_{1}(x)-v_{2}(x)|_{{\mathbb R}^{N}}
\Big(|u_{0}(x)|_{{\mathbb R}^{N}}+\frac{1}{2}|v_{1}(x)|_{{\mathbb R}^{N}}+
\frac{1}{2}|v_{2}(x)|_{{\mathbb R}^{N}}\Big).
$$
Using the Schwarz inequality, we obtain the estimate from above on
the norm $\|G_{1, m}-G_{2, m}\|_{L^{1}({\mathbb R}^{3})}$ given by
$$
M\|v_{1}-v_{2}\|_{L^{2}({\mathbb R}^{3},{\mathbb R}^{N})}\Big(\|u_{0}\|_
{L^{2}({\mathbb R}^{3},{\mathbb R}^{N})}+\frac{1}{2}\|v_{1}\|_
{L^{2}({\mathbb R}^{3},{\mathbb R}^{N})}
+\frac{1}{2}\|v_{2}\|_{L^{2}({\mathbb R}^{3},{\mathbb R}^{N})}\Big)\leq
$$
\begin{equation}
\label{g12}
M\|v_{1}-v_{2}\|_{H^{2}({\mathbb R}^{3},{\mathbb R}^{N})}
(\|u_{0}\|_{H^{2}({\mathbb R}^{3}, {\mathbb R}^{N})}+1).
\end{equation}
Thus, the upper bound on the norm
$\|u_{1, m}-u_{2, m}\|_{L^{2}({\mathbb R}^{3})}^{2}$ equals to
$$
\varepsilon^{2}\|H_{m}\|_{L^{1}({\mathbb R}^{3})}^{2}
M^{2}\|v_{1}-v_{2}\|_{H^{2}({\mathbb R}^{3}, {\mathbb R}^{N})}^{2}
\Big\{\frac{(\|u_{0}\|_{H^{2}({\mathbb R}^{3}, {\mathbb R}^{N})}+1)^{2}
R^{3-4s_{1, m}}}{2{\pi}^{2}(3-4s_{1, m})}+\frac{1}{R^{4s_{1, m}}}\Big\}.
$$
We minimize the expression above over $R\in (0,+\infty)$ via
Lemma 1.4. Hence,
$$
\|u_{1, m}-u_{2, m}\|_{L^{2}({\mathbb R}^{3})}^{2}\leq
\varepsilon^{2}\|H_{m}\|_{L^{1}({\mathbb R}^{3})}^{2}M^{2}
\|v_{1}-v_{2}\|_{H^{2}({\mathbb R}^{3}, {\mathbb R}^{N})}^{2}\times
$$
$$
(\|u_{0}\|_{H^{2}({\mathbb R}^{3}, {\mathbb R}^{N})}+1)^{\frac{8s_{1, m}}{3}}
\frac{1}{{(4s_{1, m})}^{\frac{4s_{1, m}}{3}}}
\frac{3}{(2{\pi}^{2})^{\frac{4s_{1, m}}{3}}(3-4s_{1, m})}.
$$
According to definition (\ref{es}), for $1\leq m\leq N$
$$
\|u_{1, m}-u_{2, m}\|_{L^{2}({\mathbb R}^{3})}^{2}\leq
\varepsilon^{2}\|H_{m}\|_{L^{1}({\mathbb R}^{3})}^{2}M^{2}
\|v_{1}-v_{2}\|_{H^{2}({\mathbb R}^{3}, {\mathbb R}^{N})}^{2}\times
$$
$$
(\|u_{0}\|_{H^{2}({\mathbb R}^{3}, {\mathbb R}^{N})}+1)^{\frac{8S_{1}}{3}}
\frac{1}{{(4s_{1})}^{\frac{4s_{1}}{3}}}
\frac{3}{(2{\pi}^{2})^{\frac{4s_{1}}{3}}(3-4S_{1})}.
$$
Therefore,
$$
\|u_{1}-u_{2}\|_{L^{2}({\mathbb R}^{3}, {\mathbb R}^{N})}^{2}\leq \varepsilon^{2}H^{2}
 M^{2}\|v_{1}-v_{2}\|_{H^{2}({\mathbb R}^{3}, {\mathbb R}^{N})}^{2}\times
$$
\begin{equation}
\label{u12n}
(\|u_{0}\|_{H^{2}({\mathbb R}^{3}, {\mathbb R}^{N})}+1)^{\frac{8S_{1}}{3}}
\frac{1}{{(4s_{1})}^{\frac{4s_{1}}{3}}}
\frac{3}{(2{\pi}^{2})^{\frac{4s_{1}}{3}}(3-4S_{1})}.
\end{equation}
By virtue of (\ref{aux1}) and (\ref{aux2}) with $1\leq m\leq N$, we arrive at
$$
[-\Delta+(-\Delta)^{1+s_{2, m}-s_{1, m}}](u_{1, m}(x)-u_{2, m}(x))=
$$
$$
=\varepsilon_{m}
(-\Delta)^{1-s_{1, m}}
\int_{{\mathbb R}^{3}}H_{m}(x-y)[G_{1, m}(y)-G_{2, m}(y)]dy.
$$
Let us apply the standard Fourier transform (\ref{f}) along with estimates
(\ref{fub}) and (\ref{g12}). Thus,
$$
\|\Delta (u_{1, m}-u_{2, m})\|_
{L^{2}({\mathbb R}^{3})}^{2}\leq
$$  
$$
\varepsilon^{2}\|G_{1, m}-G_{2, m}\|_{L^{1}({\mathbb R}^{3})}^{2}
\|(-\Delta)^{1-s_{1, m}}H_{m}\|_{L^{2}({\mathbb R}^{3})}^{2}\leq
$$
$$
\varepsilon^{2}M^{2}\|v_{1}-v_{2}\|_{H^{2}({\mathbb R}^{3}, {\mathbb R}^{N})}
^{2}(\|u_{0}\|_{H^{2}({\mathbb R}^{3}, {\mathbb R}^{N})}+1)^{2}
\|(-\Delta)^{1-s_{1, m}}H_{m}\|_{L^{2}({\mathbb R}^{3})}^{2}.
$$
This means that
$$
\sum_{m=1}^{N}\|\Delta (u_{1, m}-u_{2, m})\|_
{L^{2}({\mathbb R}^{3})}^{2}\leq
$$    
\begin{equation}
\label{d12}
\varepsilon^{2}M^{2}\|v_{1}-v_{2}\|_{H^{2}({\mathbb R}^{3}, {\mathbb R}^{N})}^{2}
(\|u_{0}\|_{H^{2}({\mathbb R}^{3}, {\mathbb R}^{N})}+1)^{2}Q^{2}.
\end{equation}
Bounds (\ref{u12n}) and (\ref{d12}) yield that the norm
$\|u_{1}-u_{2}\|_{H^{2}({\mathbb R}^{3}, {\mathbb R}^{N})}$ can be estimated from above by
the quantity $\varepsilon M(\|u_{0}\|_{H^{2}({\mathbb R}^{3}, {\mathbb R}^{N})}+1)\times$
\begin{equation}
\label{contr}
\Bigg\{\frac{H^{2}(\|u_{0}\|
_{H^{2}({\mathbb R}^{3}, {\mathbb R}^{N})}+1)^{\frac{8S_{1}}{3}-2}3}
{(3-4S_{1}){(2{\pi}^{2})}^{\frac{4s_{1}}{3}}}\frac{1}{(4s_{1})^{\frac{4s_{1}}{3}}}+
Q^{2}\Bigg\}^
{\frac{1}{2}}\|v_{1}-v_{2}\|_
{H^{2}({\mathbb R}^{3}, {\mathbb R}^{N})}.
\end{equation}
Note that for all the values of $\varepsilon$, which satisfy inequality
(\ref{eps}) the constant in the right side of (\ref{contr}) is less than one.
Therefore, the map $\tau_{g}: B_{\rho}\to B_{\rho}$ defined by system 
(\ref{aux}) is a strict contraction. 
Its unique fixed point $u_{p}$ is the only solution of the system of equations
(\ref{pert}) in the ball $B_{\rho}$. The resulting
$u\in H^{2}({\mathbb R}^{3}, {\mathbb R}^{N})$ given by formula (\ref{r}) solves
system (\ref{p}). Evidently, by means of (\ref{rh}), $u_{p}$ tends
to zero in the $H^{2}({\mathbb R}^{3}, {\mathbb R}^{N})$ norm as
$\varepsilon\to 0$.      \hfill\lanbox

\bigskip

We turn our attention to the proof of the second main result of the
article.

\bigskip

%%%%%%%%%%%%%%%%%%%%%%%%%%%%%%%%%%%%%%%%%%%%%%%%%%

\setcounter{section}{3}
\setcounter{equation}{0}

\centerline{\bf 3. The continuity of the resulting solution}

\bigskip

\noindent
{\it Proof of Theorem 1.5.} Obviously, for all the values of
$\varepsilon$, which satisfy (\ref{eps})
$$
u_{p,1}=\tau_{g_{1}}u_{p,1}, \quad u_{p,2}=\tau_{g_{2}}u_{p,2}.
$$
Then
$$
u_{p,1}-u_{p,2}=\tau_{g_{1}}u_{p,1}-\tau_{g_{1}}u_{p,2}+\tau_{g_{1}}u_{p,2}-
\tau_{g_{2}}u_{p,2}.
$$
Hence,
$$
\|u_{p,1}-u_{p,2}\|_{H^{2}({\mathbb R}^{3}, {\mathbb R}^{N})}\leq\|\tau_{g_{1}}u_{p,1}-\tau_{g_{1}}
u_{p,2}\|_{H^{2}({\mathbb R}^{3},  {\mathbb R}^{N})}+\|\tau_{g_{1}}u_{p,2}-\tau_{g_{2}}u_{p,2}\|_
{H^{2}({\mathbb R}^{3}, {\mathbb R}^{N})}.
$$
By virtue of inequality (\ref{contr}), we have
$$
\|\tau_{g_{1}}u_{p,1}-\tau_{g_{1}}u_{p,2}\|_{H^{2}({\mathbb R}^{3}, {\mathbb R}^{N})}\leq \varepsilon
\sigma\|u_{p,1}-u_{p,2}\|_{H^{2}({\mathbb R}^{3}, {\mathbb R}^{N})}
$$
with $\sigma$ defined in (\ref{sig}). Clearly,
$\varepsilon \sigma<1$ since the map
$\tau_{g_{1}}: B_{\rho}\to  B_{\rho}$ is a strict contraction under the given
conditions. Thus, 
\begin{equation}
\label{sigma}
(1-\varepsilon \sigma)\|u_{p,1}-u_{p,2}\|_{H^{2}({\mathbb R}^{3}, {\mathbb R}^{N})}\leq
\|\tau_{g_{1}}u_{p,2}-\tau_{g_{2}}u_{p,2}\|_{H^{2}({\mathbb R}^{3}, {\mathbb R}^{N})}.
\end{equation}
Evidently, for the fixed point we have $\tau_{g_{2}}u_{p,2}=u_{p,2}$. Let us introduce
$\gamma(x):=\tau_{g_{1}}u_{p,2}(x)$. For $1\leq m\leq N$, we arrive at
\begin{equation}
\label{12}
[(-\Delta)^{s_{1, m}}+(-\Delta)^{s_{2, m}}]\gamma_{m}(x)=\varepsilon_{m}
\int_{{\mathbb R}^{3}}
H_{m}(x-y)g_{1, m}(u_{0}(y)+u_{p,2}(y))dy,
\end{equation}
$$
[(-\Delta)^{s_{1, m}}+(-\Delta)^{s_{2, m}}]u_{p,2,m}(x)=
$$
\begin{equation}
\label{22}
\varepsilon_{m}
\int_{{\mathbb R}^{3}}H_{m}(x-y)g_{2, m}(u_{0}(y)+u_{p,2}(y))dy,
\end{equation}
where
$\displaystyle{\frac{1}{4}<s_{1, m}<\frac{3}{4}, \
s_{1, m}<s_{2, m}<1}$. We denote
$$
G_{1,2,m}(x):=g_{1,m}(u_{0}(x)+u_{p,2}(x)), \quad G_{2,2,m}(x):=g_{2,m}
(u_{0}(x)+u_{p,2}(x)).
$$
Let us apply the standard Fourier transform (\ref{f}) to both sides of
systems (\ref{12}) and (\ref{22}). This yields
$$
\widehat{\gamma_{m}}(p)=\varepsilon_{m} (2 \pi)^{\frac{3}{2}}
\frac{\widehat{H_{m}}(p)
\widehat{G_{1,2,m}}(p)}{|p|^{2s_{1, m}}+|p|^{2s_{2, m}}}, \quad
\widehat{u_{p,2,m}}(p)=\varepsilon_{m} (2 \pi)^{\frac{3}{2}}\frac{\widehat{H_{m}}
(p)\widehat{G_{2,2,m}}(p)}{|p|^{2s_{1, m}}+|p|^{2s_{2, m}}}.
$$
This means that
$$
\displaystyle{\|\gamma_{m}-u_{p,2,m}\|_{L^{2}({\mathbb R}^{3})}^{2}=}
$$
\begin{equation}
\label{xiu}  
\varepsilon_{m}^{2}{(2\pi)}^{3}
\int_{{\mathbb R}^{3}}\frac{|\widehat{H_{m}}(p)|^{2}|
\widehat{G_{1,2,m}}(p)-\widehat{G_{2,2,m}}(p)|^{2}}{[|p|^{2s_{1, m}}+|p|^{2s_{2, m}}]^{2}}
dp.
\end{equation}
We obtain the estimate from above on the right side of (\ref{xiu}) using inequality
(\ref{fub}), such that
$$
\varepsilon_{m}^{2}{(2\pi)}^{3}\Bigg[\int_{|p|\leq R}\frac{|{\widehat{H_{m}}}(p)|^{2}|\widehat{G_{1,2,m}}(p)-\widehat{G_{2,2,m}}(p)|^{2}}{|p|^{4s_{1, m}}}dp+
$$
$$  
\int_{|p|>R}\frac{|{\widehat{H_{m}}}(p)|^{2}|\widehat{G_{1,2,m}}(p)-
\widehat{G_{2,2,m}}(p)|^{2}}{|p|^{4s_{1, m}}}dp \Bigg]\leq
\varepsilon^{2}\|H_{m}\|_{L^{1}({\mathbb R}^{3})}^{2}\times
$$
$$
\Bigg\{\frac{1}
{2{\pi}^{2}}\frac{\|G_{1,2,m}-G_{2,2,m}\|_{L^{1}({\mathbb R}^{3})}^{2}R^{3-4s_{1, m}}}
{3-4s_{1, m}}+
\frac{\|G_{1,2,m}-G_{2,2,m}\|_{L^{2}({\mathbb R}^{3})}^{2}}{R^{4s_{1, m}}}\Bigg\}
$$
with $R\in (0, +\infty)$. Clearly, we can express
$$
G_{1,2,m}(x)-G_{2,2,m}(x)=\int_{0}^{1}\nabla[g_{1,m}-g_{2,m}](t(u_{0}(x)+u_{p,2}(x))).
(u_{0}(x)+u_{p,2}(x))dt.
$$
Thus,
$$
|G_{1,2,m}(x)-G_{2,2,m}(x)|\leq \|g_{1,m}-g_{2,m}\|_{C^{2}(I)}
|u_{0}(x)+u_{p,2}(x)|_{{\mathbb R}^{N}}.
$$
This implies that
$$
\|G_{1,2,m}-G_{2,2,m}\|_{L^{2}({\mathbb R}^{3})}\leq \|g_{1,m}-g_{2,m}\|_{C^{2}(I)}
\|u_{0}+u_{p,2}\|_{L^{2}({\mathbb R}^{3}, {\mathbb R}^{N})}\leq
$$
$$
\|g_{1,m}-g_{2,m}\|_{C^{2}(I)}(\|u_{0}\|_{H^{2}({\mathbb R}^{3}, {\mathbb R}^
{N})}+1).
$$
Note that for $1\leq m, j\leq N$ and $t\in [0,1]$, we have 
$$
\frac{\partial}{\partial z_{j}}(g_{1,m}-g_{2,m})(t(u_{0}(x)+u_{p,2}(x)))=
$$
$$
\int_{0}^{t}\nabla \Big[\frac{\partial}{\partial z_{j}}(g_{1,m}-g_{2,m})\Big]
(\tau(u_{0}(x)+u_{p,2}(x))).(u_{0}(x)+u_{p,2}(x))d\tau.
$$
Hence,
$$
\Big|\frac{\partial}{\partial z_{j}}(g_{1,m}-g_{2,m})(t(u_{0}(x)+u_{p,2}(x)))\Big|
\leq
$$
$$
 \sum_{n=1}^{N}\Bigg\|\frac{\partial^{2}(g_{1,m}-g_{2,m})}{\partial z_{n}
\partial z_{j}}\Bigg\|_{C(I)}|u_{0}(x)+u_{p,2}(x)|_{{\mathbb R}^{N}}.
$$
Obviously,
$$
|G_{1,2,m}(x)-G_{2,2,m}(x)|\leq \|g_{1,m}-g_{2,m}\|_{C^{2}(I)}
|u_{0}(x)+u_{p,2}(x)|_{{\mathbb R}^{N}}^{2}.
$$
Then
$$
\|G_{1,2,m}-G_{2,2,m}\|_{L^{1}({\mathbb R}^{3})}\leq
\|g_{1,m}-g_{2,m}\|_{C^{2}(I)}\|u_{0}+u_{p,2}\|_{L^{2}({\mathbb R}^{3}, {\mathbb R}^{N})}^{2}\leq
$$
\begin{equation}
\label{G1222}
\|g_{1,m}-g_{2,m}\|_{C^{2}(I)}(\|u_{0}\|_{H^{2}({\mathbb R}^{3}, {\mathbb R}^{N})}+1)^{2}.
\end{equation}
This enables us to derive the upper bound on the norm
$\|\gamma_{m}-u_{p,2,m}\|_{L^{2}({\mathbb R}^{3})}^{2}$ given by
$\varepsilon^{2}\|H_{m}\|_{L^{1}({\mathbb R}^{3})}^{2}(\|u_{0}\|_
{H^{2}({\mathbb R}^{3}, {\mathbb R}^{N})}+1)^{2} \times$
$$
\|g_{1, m}-g_{2, m}\|_{C^{2}(I)}^{2}
\Bigg[(\|u_{0}\|_{H^{2}({\mathbb R}^{3}, {\mathbb R}^{N})}+1)^{2}
\frac{R^{3-4s_{1, m}}}{2{\pi}^{2} (3-4s_{1, m})}+\frac{1}{R^{4s_{1, m}}}\Bigg].
$$
Let us minimize such expression over $R\in (0, +\infty)$ using
Lemma 1.4. This yields the estimate
$$
\|\gamma_{m}-u_{p,2,m}\|_{L^{2}({\mathbb R}^{3})}^{2}\leq
$$
$$
\varepsilon^{2}\|H_{m}\|_{L^{1}({\mathbb R}^{3})}^{2}
(\|u_{0}\|_{H^{2}({\mathbb R}^{3}, {\mathbb R}^{N})}+1)^{2+\frac{8s_{1, m}}{3}}
\frac{3\|g_{1, m}-g_{2, m}\|_{C^{2}(I)}^{2}}{(4s_{1, m})^{\frac{4s_{1, m}}{3}}(2{\pi}^{2})^{\frac{4s_{1, m}}{3}}(3-4s_{1, m})}.
$$
Therefore,
$$
\|\gamma-u_{p,2}\|_{L^{2}({\mathbb R}^{3}, {\mathbb R}^{N})}^{2}\leq
$$
$$
\varepsilon^{2}
H^{2}(\|u_{0}\|_{H^{2}({\mathbb R}^{3}, {\mathbb R}^{N})}+1)^{2+\frac{8S_{1}}{3}}
\frac{3\|g_{1}-g_{2}\|_{C^{2}(I, {\mathbb R}^{N})}^{2}}{(3-4S_{1})(2{\pi}^{2})^{\frac{4s_{1}}{3}}(4s_{1})^{\frac{4s_{1}}{3}}}.
$$
By means of (\ref{12}) and (\ref{22}) with $1\leq m\leq N$, we arrive at
$$
[-\Delta+(-\Delta)^{1+s_{2, m}-s_{1, m}}]\gamma_{m}(x)=
\varepsilon_{m}(-\Delta)^{1-s_{1, m}}
\int_{{\mathbb R}^{3}}H_{m}(x-y)G_{1,2,m}(y)dy,
$$
$$
[-\Delta+(-\Delta)^{1+s_{2, m}-s_{1, m}}]u_{p,2,m}(x)=
\varepsilon_{m}(-\Delta)^{1-s_{1, m}}
\int_{{\mathbb R}^{3}}H_{m}(x-y)G_{2,2,m}(y)dy
$$
with
$\displaystyle{\frac{1}{4}<s_{1, m}<\frac{3}{4}, \ s_{1, m}<s_{2, m}<1}$.

By virtue of the standard Fourier transform (\ref{f}) along with inequalities (\ref{fub})
and (\ref{G1222}), the norm
$\|\Delta (\gamma_{m}-u_{p,2,m})\|_{L^{2}({\mathbb R}^{3})}^{2}$
can be bounded from above by
$$
\varepsilon^{2}\|G_{1,2,m}-G_{2,2,m}\|_{L^{1}({\mathbb R}^{3})}^{2}
\|(-\Delta)^{1-s_{1, m}}H_{m}\|_
{L^{2}({\mathbb R}^{3})}^{2}\leq
$$
$$
\varepsilon^{2}\|g_{1, m}-g_{2, m}\|_{C^{2}(I)}^{2}
(\|u_{0}\|_{H^{2}({\mathbb R}^{3}, {\mathbb R}^{N})}+1)^{4}
\|(-\Delta)^{1-s_{1, m}}H_{m}\|_
{L^{2}({\mathbb R}^{3})}^{2}.
$$
Thus,
$$
\sum_{m=1}^{N}\|\Delta (\gamma_{m}-u_{p,2,m})
\|_{L^{2}({\mathbb R}^{3})}^{2}\leq
$$
$$
\varepsilon^{2}
\|g_{1}-g_{2}\|_{C^{2}(I,{\mathbb R}^{N})}^{2}
(\|u_{0}\|_{H^{2}({\mathbb R}^{3}, {\mathbb R}^{N})}+1)^{4}Q^{2}.
$$
This means that
$$
\|\gamma-u_{p,2}\|_
{H^{2}({\mathbb R}^{3}, {\mathbb R}^{N})}\leq \varepsilon
\|g_{1}-g_{2}\|_{C^{2}(I, {\mathbb R}^{N})}\times
$$
$$
(\|u_{0}\|_{H^{2}({\mathbb R}^{3}, {\mathbb R}^{N})}+1)^{2}\Bigg[
\frac{H^{2}(\|u_{0}\|_{H^{2}({\mathbb R}^{3},{\mathbb R}^{N})}+1)^{\frac{8S_{1}}{3}-2}3}
{(3-4S_{1})
(2{\pi}^{2})^{\frac{4s_{1}}{3}}(4s_{1})^{\frac{4s_{1}}{3}}}+Q^{2}\Bigg]^
\frac{1}{2}.
$$
Let us recall estimate (\ref{sigma}). The norm
$\|u_{p,1}-u_{p,2}\|_{H^{2}({\mathbb R}^{3}, {\mathbb R}^{N})}$ can be bounded from above
by
$$
\frac{\varepsilon}{1-\varepsilon \sigma}
(\|u_{0}\|_{H^{2}({\mathbb R}^{3}, {\mathbb R}^{N})}+1)^{2}\times
$$
$$
\Bigg[\frac{{H}^{2}
(\|u_{0}\|_{H^{2}({\mathbb R}^{3}, {\mathbb R}^{N})}+1)^{\frac{8S_{1}}{3}-2}3}{(3-4S_{1})
(2{\pi}^{2})^{\frac{4s_{1}}{3}}(4s_{1})^{\frac{4s_{1}}{3}}}+
Q^{2}\Bigg]^{\frac{1}{2}}\|g_{1}-g_{2}\|_{C^{2}(I, {\mathbb R}^{N})}.
$$
By means of equalities (\ref{sig}) and (\ref{jres}) we complete the proof of our
theorem. \hfill\lanbox

\bigskip

%%%%%%%%%%%%%%%%%%%%%%%%%%%%%%%%%%%%%%%%%%%%%%%%%%

\setcounter{section}{4}
\setcounter{equation}{0}

\centerline{\bf 4. Auxiliary results}

\bigskip
\bigskip

Let us derive the solvability conditions for the linear
Poisson type equation with a square integrable right side in the case
of the double scale anomalous diffusion
\begin{equation}
\label{lp}
[{(-\Delta)}^{s_{1}}+{(-\Delta)}^{s_{2}}]\phi(x)=f(x), \quad x\in {\mathbb R}^{3},
\quad 0<s_{1}<s_{2}<1.
\end{equation}
The inner product is given by
\begin{equation}
\label{ip}
(f(x), g(x))_{L^{2}({\mathbb R}^{3})}:=\int_{{\mathbb R}^{3}}f(x){\bar g}(x)dx,
\end{equation}
with a slight abuse of notations when the functions involved in (\ref{ip})
are not square integrable, like for instance the one present in
orthogonality relation (\ref{oc}) below. Clearly, for
$f(x)\in L^{1}({\mathbb R}^{3})$ and $g(x)\in L^{\infty}({\mathbb R}^{3})$, the
integral in the right side of formula (\ref{ip}) makes sense.

The technical proposition below was established in the preceding article
~\cite{V24} via the standard Fourier transform (\ref{f}). We provide
the argument below for the convenience of the readers.

\bigskip

\noindent
{\bf Lemma 4.1.} {\it  Let  $0<s_{1}<s_{2}<1, \ f(x): {\mathbb R}^{3}\to
{\mathbb R}, \ f(x)\in L^{2} ({\mathbb R}^{3})$.

\medskip

\noindent
a) Suppose $\displaystyle{s_{1}\in \Big(0, \frac{3}{4}\Big)}$ and in addition
$f(x)\in L^{1} ({\mathbb R}^{3})$.
Then equation (\ref{lp}) has a unique solution
$\phi(x)\in H^{2s_{2}}({\mathbb R}^{3})$.

\medskip

\noindent
b) Suppose $\displaystyle{s_{1}\in \Big[\frac{3}{4} , 1\Big)}$ and additionally
$xf(x)\in L^{1} ({\mathbb R}^{3})$. Then problem (\ref{lp}) possesses
a unique solution $\phi(x)\in H^{2s_{2}}({\mathbb R}^{3})$ if and only if
the orthogonality condition
\begin{equation}
\label{oc}
(f(x), 1)_{L^{2} ({\mathbb R}^{3})}=0
\end{equation}
is valid.}

\bigskip

\noindent
{\it Proof.} It can be easily verified that if
$\phi(x)\in L^{2} ({\mathbb R}^{3})$ is a solution of equation (\ref{lp}) with
a square integrable right side, it will belong to
$H^{2s_{2}} ({\mathbb R}^{3})$ as well. To demonstrate that, we apply the
standard Fourier transform (\ref{f}) to both sides of (\ref{lp}). This yields
$$
(|p|^{2s_{1}}+|p|^{2s_{2}})\widehat{\phi}(p)=\widehat{f}(p)\in
L^{2} ({\mathbb R}^{3}).
$$
Thus,
$$
\int_{{\mathbb R}^{3}}[|p|^{2s_{1}}+|p|^{2s_{2}}]^{2}|\widehat{\phi}(p)|^{2}dp<\infty.
$$
We have the trivial identity
$$
\|(-\Delta)^{s_{2}}\phi\|_{L^{2} ({\mathbb R}^{3})}^{2}=\int_{{\mathbb R}^{3}}|p|^{4s_{2}}
|\widehat{\phi}(p)|^{2}dp<\infty.
$$
This means that $(-\Delta)^{s_{2}}\phi\in L^{2} ({\mathbb R}^{3})$. We use
the definition of the norm (\ref{n}). Therefore, 
$\phi(x)\in H^{2s_{2}} ({\mathbb R}^{3})$ as well.

Let us suppose that equation (\ref{lp}) has admits two solutions
$\phi_{1, 2}(x)\in H^{2s_{2}} ({\mathbb R}^{3})$. Then their difference
$w(x):=\phi_{1}(x)-\phi_{2}(x)\in H^{2s_{2}} ({\mathbb R}^{3})$ is a solution
of the homogeneous problem
$$
[(-\Delta)^{s_{1}}+(-\Delta)^{s_{2}}]w=0.
$$
Evidently, the operator
$$
(-\Delta)^{s_{1}}+(-\Delta)^{s_{2}}: H^{2s_{2}} ({\mathbb R}^{3})\to
L^{2} ({\mathbb R}^{3})
$$
does not have any nontrivial zero modes. Thus, $w(x)$ vanishes in
${\mathbb R}^{3}$, which yields the uniqueness of solutions for problem
(\ref{lp}).

We apply the standard Fourier transform (\ref{f}) to both sides of
equation (\ref{lp}) and arrive at
\begin{equation}
\label{fihpfp}
\widehat{\phi}(p)=\frac{\widehat{f}(p)}{|p|^{2s_{1}}+|p|^{2s_{2}}}\chi_{\{|p|\leq 1\}}+
\frac{\widehat{f}(p)}{|p|^{2s_{1}}+|p|^{2s_{2}}}\chi_{\{|p|>1\}}.
\end{equation}
In formula (\ref{fihpfp}) and further down $\chi_{A}$ will designate the
characteristic function of a set $A\subseteq {\mathbb R}^{3}$.

Clearly, the second term in the right side of (\ref{fihpfp}) can be bounded
from above in the absolute value by
$\displaystyle{\frac{|\widehat{f}(p)|}{2}\in L^{2}({\mathbb R}^{3})}$ via 
the one of our assumptions.

The first term in the right side of (\ref{fihpfp}) can be trivially estimated
from above in the absolute value using (\ref{fub}) by
\begin{equation}
\label{fxl1p2s2}  
\frac{\|f(x)\|_{L^{1}({\mathbb R}^{3})}}{(2\pi)^{\frac{3}{2}}|p|^{2s_{1}}}
\chi_{\{|p|\leq 1\}}.
\end{equation}
It can be trivially checked that quantity (\ref{fxl1p2s2}) is square
integrable in our space of three dimensions when
$\displaystyle{s_{1}\in \Big(0, \frac{3}{4}\Big)}$.

Let us complete the argument by considering the case when
$\displaystyle{s_{1}\in \Bigg[\frac{3}{4}, 1\Bigg)}$. Note that
\begin{equation}
\label{fhp0i}
\widehat{f}(p)=\widehat{f}(0)+\int_{0}^{|p|}
\frac{\partial \widehat{f}(q, \sigma)}{\partial q}dq.
\end{equation}
Here and below $\sigma$ will stand for the angle variables on the sphere. Hence,
the first term in the right side of (\ref{fihpfp}) is given by
\begin{equation}
\label{fh0ps1s2}
\frac{\widehat{f}(0)}{|p|^{2s_{1}}+|p|^{2s_{2}}}\chi_{\{|p|\leq 1\}}+
\frac{\int_{0}^{|p|}\frac{\partial \widehat{f}(q, \sigma)}{\partial q}dq}
{|p|^{2s_{1}}+|p|^{2s_{2}}}\chi_{\{|p|\leq 1\}}.    
\end{equation}
We use the definition of the standard Fourier transform (\ref{f}) to derive
that
\begin{equation}
\label{ftdub}
\Bigg|\frac{\partial \widehat{f}(p)}{\partial |p|}\Bigg|\leq
\frac{\|xf(x)\|_{L^{1}({\mathbb R}^{3})}}{(2\pi)^{\frac{3}{2}}}.
\end{equation}
Therefore,
$$
\Bigg|\frac{\int_{0}^{|p|}\frac{\partial \widehat{f}(q, \sigma)}{\partial q}dq}
{|p|^{2s_{1}}+|p|^{2s_{2}}}\chi_{\{|p|\leq 1\}}\Bigg|\leq     
\frac{\|xf(x)\|_{L^{1}({\mathbb R}^{3})}}{(2\pi)^{\frac{3}{2}}}|p|^{1-2s_{1}}
\chi_{\{|p|\leq 1\}}\in L^{2}({\mathbb R}^{3}).
$$
Let us analyze the remaining term
\begin{equation}
\label{fh0ps12}  
\frac{\widehat{f}(0)}{|p|^{2s_{1}}+|p|^{2s_{2}}}\chi_{\{|p|\leq 1\}}.
\end{equation}
Clearly, (\ref{fh0ps12}) is square integrable in the space of three dimensions
if and only if $\widehat{f}(0)$ is trivial. This is equivalent to orthogonality
relation (\ref{oc}).  \hfill\lanbox.

\bigskip

The result of the technical lemma above deals with the solvability of problem
(\ref{lp}) in $H^{2s_{2}}({\mathbb R}^{3})$
for all the values of the powers of the fractional Laplace operators
$\displaystyle{0<s_{1}<s_{2}<1}$, so that no orthogonality conditions are
needed for the right side $f(x)$ for
$\displaystyle{s_{1}\in \Bigg(0, \frac{3}{4}\Bigg)}$. But the argument of
our proof is based on the single orthogonality condition (\ref{oc}) when
$\displaystyle{s_{1}\in \Bigg[\frac{3}{4}, 1\Bigg)}$.
This is analogous to the situation when the
Poisson type equation is considered with a single fractional Laplacian in
${\mathbb R}^{3}$ (see Theorem 1.1 of ~\cite{VV19}). The
solvability of the equation similar to (\ref{lp})
involving a scalar potential was discussed in ~\cite{EV22} and ~\cite{V24}.

\bigskip

%%%%%%%%%%%%%%%%%%%%%%%%%%%%%%%%%%%%%%%%%%%%%%%%%%%%%%%%%%%%%%%%

\section*{Acknowledgements}

Vitali Vougalter is grateful to Israel Michael Sigal for the partial support
by the NSERC grant NA 7901. Vitaly Volpert has been supported by the RUDN
University Strategic Academic Leadership Program.

\bigskip

%%%%%%%%%%%%%%%%%%%%%%%%%%%%%%%%%%%%%%%%%%%%%%%%%%%%%%%%%%%%%%%%

\end{document}